\documentclass{article}

\usepackage[a4paper,top=2.54cm,bottom=2.54cm,left=3.17cm,right=3.17cm,%
            includehead,includefoot]{geometry}

\usepackage{amsmath,amssymb,amsfonts,amsthm}
\usepackage{graphicx}
\usepackage{subfigure}
\usepackage{float}
\usepackage{xcolor}
\usepackage[numbers,square,sort&compress]{natbib}
\usepackage{hyperref}
  \hypersetup{colorlinks,citecolor=blue,linkcolor=blue,breaklinks=true}
\usepackage{booktabs}
\usepackage{colortbl}
\usepackage{caption}
\usepackage{enumitem}
\usepackage{epstopdf}
\usepackage{algorithm}
\usepackage{algpseudocode}
\usepackage{array}
\usepackage{bbding}
\usepackage{fancyhdr} 
\usepackage{fancyvrb} 
\usepackage{longtable}
\usepackage{listings}
\usepackage{rotating,rotfloat} 
\usepackage{yhmath} 
\usepackage{mathrsfs}
\usepackage{multirow}
\usepackage{makecell}
\usepackage{graphicx}
\usepackage{fancyhdr}
\usepackage{setspace}
\allowdisplaybreaks

\newtheorem{theorem}{Theorem}[section]
\newtheorem{assumption}{Assumption}
\newtheorem{condition}{Condition}

\newtheorem{remark}{Remark}
\newcommand{\keywords}[1]{\par\noindent\textbf{Keywords:} #1}

\newcommand{\bbm}{\begin{bmatrix}}
\newcommand{\ebm}{\end{bmatrix}}

\begin{document}

\title{\uppercase{Group Efficient Randomized-Adaptive Designs With Delayed and Missing Responses}}

\author{
  Guijing Zhang
  \and
  Lixin Zhang
    \thanks{School of Statistics and Data Science, Zhejiang Gongshang University, Zhejiang, China}
}

\date{}

\maketitle

\begin{abstract}
	Response-adaptive randomization designs have attracted much attention in clinical trials. This paper proposes a new class of response-adaptive design which built on the efficient randomized-adaptive design (ERADE) proposed by Hu, Zhang, and He. The original ERADE uses a discrete allocation probability function and leverages the stopping time theory of stochastic processes to establish asymptotic results. It has been proven that the design can reach the Cramér-Rao lower bound for any target allocation proportion. This study further expands the original framework by replacing the traditional case-by-case sequential enrollment with group recruitment over fixed time intervals(weekly, biweekly, or monthly), and dynamically updates the allocation probabilities based on the cumulative response information within each group. Meanwhile, to better fit practical application scenarios, we further explicitly consider the situations of randomly missing data and response delay. Theoretical analysis shows that the new design retains all the main asymptotic properties of the original ERADE, and still performs well under the conditions of response delay and missing data. Finally, through simulation studies and the redesign of a real-world clinical trial, the effectiveness and practicality of the proposed method are verified.
\end{abstract}

\section{Preliminaries}\label{sec1}

In clinical trials, patients typically enroll sequentially over a specified period. Researchers can adjust allocation probability using accumulated information, such as treatment assignments and outcomes, to achieve a predefined target allocation proportion. This approach is termed a randomized-adaptive design. The three core components of a randomized-adaptive design are allocation proportion, efficiency, and variability\cite{hu2009efficient}.

The concept of randomized-adaptive design can be traced back to Thompson \cite{thompson1933likelihood} and Robbins \cite{robbins1952some}, whose core idea was to allocate more patients to superior treatments. Subsequently, numerous scholars expanded on this randomized-adaptive framework. Wei and Durham \cite{wei1978randomized} proposed the randomized play-the-winner rule (RPW), and Ivanova \cite{ivanova2003play} further developed the drop-the-loser rule (DL) based on this foundation. Compared to the RPW strategy, DL exhibits lower variability. Hu and Rosenberger \cite{hu2003optimality} proposed that in randonzied-adaptive designs, the variance of the actual allocation proportion is commonly used to quantify the variability of the design.

The optimal randomized-adaptive design is a response-driven randomized-adaptive design implemented to achieve a specific optimal allocation proportion. This proportion typically depends on a defined proportional function of unknown parameters. For instance, Rosenberger et al. \cite{rosenberger2001optimal} proposed an allocation proportion that minimizes the number of treatment failures. Notably, the allocation proportion of the randomized urn model generally does not qualify as the optimal allocation in the context of statistical hypothesis testing. Additionally, Eisele \cite{eisele1994doubly}, Eisele and Woodroofe \cite{eisele1995central} introduced the doubly adaptive biased coin design (DBCD) and demonstrated that, under stringent conditions, dynamically adjusting allocation probabilities ensures convergence of the actual allocation proportion to the desired allocation proportion. Hu and Zhang \cite{hu2004asymptotic} further developed a set of DBCD allocation functions and derived their asymptotic properties under more generalized constraints. Hu, Rosenberger and Zhang \cite{hu2006asymptotically} proved that while the variance of DBCD’s actual allocation proportion does not attain the Cramér-Rao loser bound, the discrepancy is negligible. Hu, Zhang and He \cite{hu2009efficient} proposed the efficient randomized-adaptive design (ERADE), whose allocation proportion is unbiased and achieves the Cramér-Rao bound, establishing it as the theoretically optimal randomized-adaptive design.

In actual clinical trials, missing and delayed responses are common problems that cannot be completely avoided. However, in existing research on response-adaptive randomization design, most methods directly ignore missing data or delayed responses, and only update allocation probabilities based on the observed complete data. Although this approach simplifies the analytical process, it may introduce systematic bias, thereby undermining the reliability and validity of trial results. In response to this problem, some scholars have conducted relevant theoretical explorations. Bai, Hu, and Rosenberger \cite{bai2002asymptotic} conducted a systematic theoretical analysis of the response delay mechanism in the urn model, laying an important foundation for understanding the impact of delayed data on adaptive designs. Subsequently, Hu et al. \cite{hu2008doubly} further studied the impact of response delays on DBCD and derived several asymptotic properties of this design under delayed response. Zhai et al. \cite{zhai2024group} verified the effectiveness of the group DBCD in the presence of delayed and missing responses, and the results showed that its performance was basically consistent with the traditional DBCD design.

In this paper, to address the operational challenges posed by frequent updates of allocation probabilities in response-adaptive designs, we proposes a novel design based on ERADE: group efficient randomized-adaptive design (Group ERADE). This design recruits subjects in groups at fixed time intervals (e.g., weekly, biweekly, or monthly), replacing the traditional case-by-case aequential enrollment. It dynamically updates the treatment allocation probabilities for the next group based on the accumulated response information within each group. On this basis, to handle practical issues such as response delays and missing data commonly found in clinical trials, this paper further proposes an improved version of Group ERADE, enabling it to operate effectively under conditions closer to real-world application scenarios. Theoretically, by introducing the martingale process and its associated stopping times, we prove that both Group ERADE and the improved version considering delayed and missing data retain the basic theoretical properties of ERADE. Specifically, under mild conditions, the allocation proportions of both designs exhibit asymptotic normality and strong consistency.

The paper is structured as follows. Section 2 starts with the framework of the proposed Group ERADE with delayed and missing responses and discusses its asymptotic properties. In Section 3, we apply our design to a clinical trial evaluating the efficacy and safety of a drug and examine its performance under finite samples through extensive simulation studies. Finally, Section 4 deals with some conclusions. Supplementary theoretical details and additional simulation results are given in the appendix.

\section{Group ERADE with Delayed and Missing Responses}\label{sec2}

In this section, we propose an important generalization of ERADE, namely Group ERADE with delayed and missing responses. Unlike the traditional ERADE, this procedure allows sequential group-wise patient enrollment and response observation, while fully accounting for two common complex scenarios in practical clinical trials: treatment response delay (for example, it may take several days or weeks to obtain outcomes) and missing response data caused by early participant withdrawal or other reasons. We first elaborate on the basic framework and implementation procedure of this new design in detail, then systematically discuss its corresponding asymptotic theoretical properties, including the strong consistency and asymptotic normality of the allocation proportion, to demonstrate that this generalized design can maintain robust statistical performance under more realistic clinical scenarios.

\subsection{Framework}

We consider two-treatment clinical trials. Suppose that the sizes of the groups are assumed to be $\chi_1$, $\chi_2$, … , which are random integers with $\chi_j\geq1$. Let $\tau_j=\sum_{m=1}^{j}\chi_m$ be the number of patients assigned to the first $j$ groups. Regarding the group sizes $\chi_j$, we make the following assumption.

\begin{assumption}
	For some $\varphi>0$,
	$$\chi_j = o\left(\tau_j^{1/2-\varphi}\right)\quad \text{a.s.}$$
	This assumption holds if the $(2+\varepsilon_0)$ moments $\mathbb{E}\left[\chi_j^{2+\varepsilon_0}\right]$ are bounded.
	Specifically, for any $0<\varphi<\frac{1}{2\left(2+\varepsilon_0\right)}$, we have
	$$\sum_{j} \mathbb{P}\left(\chi_j\geq\varepsilon j^{1/2-\varphi}\right)\leq\sum_{j}\frac{\mathbb{E}\left[\chi_j^{2+\varepsilon_0}\right]}{\left(\varepsilon j^{1/2-\varphi}\right)^{2+\varepsilon_0}}=\sum_{j} Cj^{-\left(1/2-\varphi \right)\left(2+\varepsilon_0\right)}<\infty$$
	By the Borel–Cantelli lemma and the fact that $\tau_j\ge j$, it follows that
	$$\chi_j = o\left(\tau_j^{1/2-\varphi}\right)\quad \text{a.s.}$$
\end{assumption}

Let $\mathbf{X}_{k}^{\left(j\right)}=\left(X_{1,k}^{\left(j\right)},X_{2,k}^{\left(j\right)},\cdots,X_{\chi_j,k}^{\left(j\right)}\right)$ be the allocation vector for the $j$th group of patients (i.e., $X_k^{\left(j\right)} = 1$ if the $j$th group of patients is assigned to treatment $k$, and $0$ otherwise), and let $N_{\tau_j,k}=\sum_{m=1}^{j} \sum_{i=1}^{\chi_m} X_{i,k}^{\left(m\right)}$ denote the number of patients assigned to treatment $k$ up to the $j$th group. Let $N_{m,k}$ be the number of patients assigned to treatment $k$ in the first assignments, i.e., when $m=\tau_j+l\le \tau_{j+1}$, $N_{m,k}=N_{\tau_j,k}+\sum_{i=1}^l X_{i,k}^{j+1}$. We denote $\xi_{i,k}^{\left(j\right)}$ as the response of the $i$th patient in the $j$th group to treatment $k$. We assume that the responses of each patient to each treatment are independent and identically distributed, following a probability distribution indexed by the parameter $\theta_k\in \mathcal{R}^d\left(k=1,2\right)$. 

Next, we introduce some additional notation related to delayed and missing responses. Let $t_i^{\left(j\right)}$ represent the enrollment time of the $i$th patient in the $j$th group. Within the same group, the enrollment times of all patients are equal, i.e., $t^{\left(j\right)}=t_{1}^{\left(j\right)}=\cdots=t_{\chi_j}^{\left(j\right)}$, and the condition $t^{\left(j\right)}<t^{\left(j+1\right)}$ holds, thus forming a strictly increasing time sequence (for example, patients are enrolled at weekly, biweekly or monthly intervals). For the $i$th patient in the $j$th group, the response delay time under treatment $k$ is a non-negative random variable, denoted as $\eta_{i,k}^{\left(j\right)}$, where $k=1,2$. To describe whether the response is observed within a specific time window, the following indicator function is introduced: $M_{i,k}\left(j,l\right)=I\left(t^{\left(j\right)}+\eta_{i,k}^{\left(j\right)}\in \left[t^{\left(j+l\right)},t^{\left(j+l+1\right)}\right)\right)$. A value of $1$ for this function indicates that the response $\xi_{i,k}^{\left(j\right)}$ of the $i$th patient in the $j$th group under treatment $k$ is observed between the treatment allocation of the $j+l$th group and the $j+l+1$th group; otherwise, the value is $0$. In addition, to handle the missing data problem, we define another indicator function $I_{mis,k}^{\left(j\right)}\left(i\right)$ to mark whether the response $\xi_{i,k}^{\left(j\right)}$ is missing. This variable follows a Bernoulli distribution with parameter $\beta_k \left(0\leq\beta_k<1\right)$ and is independent of the response value $\xi_{i,k}^{\left(j\right)}$. Before the treatment allocation of the $j+1$th group, the cumulative number of valid responses to treatment $k$ and the cumulative sum of their response values are defined as:
$$T_{\tau_j,k}=\sum_{g=1}^{j}\sum_{l=0}^{j-g}\sum_{i=1}^{\chi_g}M_{i,k}\left(g,l\right)\left(1-I_{mis,k}^{\left(g\right)}\left(i\right)\right)X_{i,k}^{\left(g\right)}$$
$$S_{\tau_j,k}=\sum_{g=1}^{j}\sum_{l=0}^{j-g}\sum_{i=1}^{\chi_g}M_{i,k}\left(g,l\right)\left(1-I_{mis,k}^{\left(g\right)}\left(i\right)\right)X_{i,k}^{\left(g\right)}\xi_{i,k}^{\left(g\right)}$$

In clinical trials, one of the core objectives is to ensure that the actual allocation proportion $N_{\tau_j,k}/{\tau_j}$ of each treatment group converges in probability to the pre-specified target allocation proportion. In this paper, we assume that the target allocation proportion can be expressed as a function of the parameter $\boldsymbol\Theta=\left(\theta_1,\theta_2\right)$, where $\theta_k$ reflects the response characteristics of patients under treatment $k$.\cite{rosenberger2001optimal} In practical implementation, researchers usually need to provide an initial estimate of $\boldsymbol\Theta$ in advance, denoted as $\boldsymbol\Theta_0$. This initial value can be reasonably inferred from existing literature or previous trial results, or derived from historical data of similar populations. In many practical applications, the distribution of the response variable typically belongs to the exponential family. Under this framework, we further assume that the parameter $\theta_k$ represents the expected value of the response observed by patients receiving treatment $k$, that is, $\theta_k=\mathbb{E}\left[\boldsymbol{\xi}_{1,k}^{\left(1\right)}\right]$. We define $\theta_{j,k}^{dm}$ as the response characteristics of the $j$th group of patients under treatment $k$ in the presence of delayed and missing responses. To effectively estimate $\theta_{j,k}^{dm}$, this paper adopts the following modified sample mean estimator:
$$\widehat{\theta}_{j,k}^{dm}=\dfrac{S_{\tau_j,k}+\theta_{0,k}^{dm}}{T_{\tau_j,k}+1}$$

The Group ERADE procedure with the delayed and missing responses is defined in the following:

Initial Step: To begin, patients in the first two groups were assigned to different treatment regimens using a restricted randomization method to ensure balanced treatment allocation.

Estimation Step: Assume that the first $j$ groups of patients have been assigned in the trial. Let $\widehat{\boldsymbol\Theta}_j^{dm}=\left(\widehat{\theta}_{j,1}^{dm},\widehat{\theta}_{j,2}^{dm}\right)$ denote the estimator of $\boldsymbol\Theta_j^{dm}=\left(\theta_{j,1}^{dm},\theta_{j,2}^{dm}\right)$ derived from these $\tau_j$ responses, and let $\widehat{\rho}_j=\rho\left(\widehat{\boldsymbol\Theta}_j^{dm}\right)$, where $\rho\left(\cdotp\right): \mathcal{R}^{d\times2} \to \left(0,1\right)$ is an allocation proportion.

Allocation Step: After the ﬁrst $j$ groups have been assigned to treatments, the $j+1$th group will be assigned to treatment 1 with probability $p_{j+1}^{dm}$ and to treatment $2$ with probability $1-p_{j+1}^{dm}$. The probability $p_{j+1}^{dm}$ may depend on both the treatment regimens assigned to the first $j$ groups and the observed responses of the first $\tau_j$ patients. Then, we show that Group ERADE procedure with the delayed and missing responses assigns the patients in the $j+1$th group to treatment $1$ with the following probability:
$$p_{j+1}^{dm}=
\begin{cases}
	\alpha \widehat{\rho}_j^{dm}, & \text{if } N_{\tau_j,1}/{\tau_j} > \widehat{\rho}_j^{dm},\\
	\widehat{\rho}_j^{dm}, & \text{if } N_{\tau_j,1}/{\tau_j} = \widehat{\rho}_j^{dm},\\
	1 - \alpha\left(1 - \widehat{\rho}_j^{dm}\right), & \text{if } N_{\tau_j,1}/{\tau_j} < \widehat{\rho}_j^{dm},
\end{cases}$$

where $0\leq\alpha<1$ is a constant reﬂecting the degree of randomization.

Repeated Step: The estimation step and the allocation step are then repeated iteratively until the end of the clinical trial.

When $\chi_j=1$(i.e., each group consists of only one patient) and there is no delayed or missing response, the proposed design degenerates to ERADE.\cite{hu2009efficient}

\begin{remark}
	It is easy to see that $\alpha$ is related to the randomness of the design. When $\alpha$ is small, Group ERADE is more deterministic and may have less variability. In the allocation formula, if $\alpha=2/3$ and $\widehat{\rho}_j=1/2$, the design corresponds to Efron’s biased coin design (BCD) \cite{efron1971forcing}. Efron first proposed the BCD, clarifying the bias rule for $\alpha$ values under experiments imbalance and laying a theoretical framework for randomized allocation. Wei \cite{wei1978adaptive} proposed the generalized biased coin design, derived the relationship between $\alpha$ and sample size imbalance, and facilitated the transition of $\alpha$ from fixed to dynamically adjusted values. Burman  \cite{burman1996sequential} introduced the expected p-value deficiency to assess the performance of a given design; building on this work, selecting $\alpha$ between 0.4 and 0.7 is reasonable. Hu, Zhang and He \cite{hu2009efficient} verified this conclusion through simulation studies.
\end{remark}

\subsection{Asymptotic Properties}

In this section, the asymptotic properties of Group ERADE with delayed and missing responses proposed in Section 2.1 will be presented. First, we state some conditions and assumptions.

\begin{condition}
	In the Bahadur-type representation: 
	$$\widehat{\theta}_{j,k}=N_{\tau_j,k}^{-1} \sum_{m=1}^{j}\sum_{i=1}^{\chi_m}X_{i,k}^{\left(m\right)}{\xi}_{i,k}^{\left(m\right)}+o\left(N_{\tau_j,k}^{-1/2}\right)\quad \text{as } \tau_j \to \infty$$
	The expected value of the response variable is $\theta_k=\mathbb{E}\left[\boldsymbol{\xi}_{1,k}^{\left(1\right)}\right]$, and its variance is finite, i.e., $\mathbb{E}\left\|\boldsymbol{\xi}_{1,k}^{\left(1\right)}\right\|^2<\infty$ for $k=1,2$.
\end{condition}

\begin{condition}
	The proportion function $\rho\left(\cdotp\right): \mathcal{R}^{d\times2} \to \left(0,1\right)$, defined by $\mathbf{y} \to \rho\left(\mathbf{y}\right)$ is a continuous twice differentiable at $\boldsymbol{\Theta}$ with $\rho\left(\boldsymbol{\Theta}\right)=\nu$.
\end{condition}

\noindent Let
\begin{align*}
	\mathcal{G}_{\tau_j}=\sigma\Bigl(&X_{1,1}^{(1)},\ldots,X_{\chi_j,1}^{(j)},\xi_1^{(1)},\ldots,\xi_{\chi_j}^{(j)},I_{mis,k}^{(1)}(1),\ldots,I_{mis,k}^{(j)}(\chi_j),\\
	&M_{1,k}(1,l),\ldots,M_{\chi_j,k}(j,l)\Bigr)
\end{align*}
be the sigma field generated by the previous $j$ stages. Then under $\mathcal{G}_{\tau_{\left(j-1\right)}}$, $\xi_{\chi_j}^{\left(j\right)}$, $I_{mis,k}^{\left(j\right)}\left(\chi_j\right)$ and $M_{\chi_j,k}\left(j,l\right)$ are independent. In the absence of missing data, that is, when $I_{mis,k}^{\left(j\right)}\left(i\right)=0$, although the patient's treatment response may have a certain delay, it will eventually be observed within a certain time window. At this time, the indicator function $M_{i,k}\left(j,l\right)$ takes the value $1$ for exactly one of all possible delay windows $l$, that is, it satisfies: $\sum_{l}M_{i,k}\left(j,l\right)=1$. This property indicates that as long as the response is not missing, regardless of how long the delay time is, it will eventually be observed in some finite time interval. To characterize the probability structure of this delay mechanism, we define the following conditional expectation:
$$\zeta_{l,k}^{\left(j\right)}\left(i\right)=\mathbb{E}\left(M_{i,k}\left(j,l\right) \mid I_{mis,k}^{\left(j\right)}\left(i\right)=0\right)$$
This quantity represents the probability that the response is exactly observed in the $l$th time window, conditional on the response not being missing. It follows directly from the above property that the sum of these probabilities over all possible windows equals $1$, that is: $\sum_{l}\zeta_{l,k}^{\left(j\right)}\left(i\right)=1$. Based on the above setup, we further assume that on the time scale of the entire trial period, the response delay time will not be excessively long.

\begin{assumption}\label{Ass2}
	For some $\phi>0$, then
	$$\sum_{q=l}^{\infty}\zeta_{q,k}^{\left(j\right)}\left(i\right)= \mathbb{P}\left(\eta_{i,k}^{\left(j\right)}>t_{\chi_{j+l}}^{\left(j+l\right)}-t_{\chi_j}^{\left(j\right)} \mid I_{mis,k}^{\left(j\right)}\left(i\right)=0\right)=o\left(\tau_l^{-\phi}\right)$$
	Notice $t^{\left(j\right)}=t_{1}^{\left(j\right)}=\cdots=t_{\chi_j}^{\left(j\right)}<t^{\left(j+1\right)}$, hence
	$$\sum_{q=l}^{\infty}\zeta_{q,k}^{\left(j\right)}\left(i\right)= \mathbb{P}\left(\eta_{i,k}^{\left(j\right)}>t^{\left(j+l\right)}-t^{\left(j\right)} \mid I_{mis,k}^{\left(j\right)}\left(i\right)=0\right)=o\left(l^{-\phi}\right)$$
\end{assumption}

We now start the main asymptotic properties of the Group ERADE with delayed and missing responses. For convenience of notation, define ${V}_k=\mathrm{Var}\left(\boldsymbol{\xi}_{1,k}^{\left(1\right)}\right)$ as the marginal response variance of treatment $k$ for $k=1,2$, and further introduce the following variance matrix:
$$\mathbf{V}^{mis}=\mathrm{diag}\left(\frac{1}{\left(1-\beta_1\right)\nu}{V}_1,\frac{1}{\left(1-\beta_2\right)\left(1-\nu\right)}{V}_2\right)$$
Correspondingly, the asymptotic variance is:
$${\sigma^2}^{mis} = \left(\left.\frac{\partial \rho}{\partial \mathbf{y}}\right|_{\boldsymbol{\Theta}}\right)'\mathbf{V}^{mis} \left.\frac{\partial \rho}{\partial \mathbf{y}}\right|_{\boldsymbol{\Theta}}$$

\begin{theorem}\label{thm1}
	Under Conditions 1 and 2, and when Assumption 1 and 2 are satisfied, and let $\varphi'=\left(2+\varepsilon_0\right)/\left(2+\varepsilon_0+\left(1+\varepsilon_0\right)\phi\right)$. Then,
	$$\widehat{\boldsymbol{\Theta}}_j^{dm}-\boldsymbol{\Theta}=O\left(\sqrt{\frac{\log \log \tau_j}{\tau_j}}\right)+o\left(\tau_j^{\varphi'-1}\left(\log \tau_j\right)^2\right)\quad \text{a.s.}$$
	Furthermore, if $\phi>\left(2+\varepsilon_0\right)/\left(1+\varepsilon_0\right)$, we have, as $\tau_j \to \infty$
	$$\left|N_{\tau_j,1}-\tau_j\widehat{\rho}_j^{dm}\right|=O\left(\sqrt{\tau_j\log \log\tau_j}\right)$$
	$$\sqrt{\tau_j}\left(N_{\tau_j,1}/\tau_j-\nu\right) \xrightarrow{D}\mathrm{N}\left(0,{\sigma^2}^{mis}\right)$$
\end{theorem}

\begin{theorem}
	Under Conditions 1 and 2, and when Assumption 1 and 2 are satisfied, let $I_k$ denote the Fisher information matrix for parameter $\theta_k\left(k=1,2\right)$. if $\mathrm{Var}\left(\boldsymbol{\xi}_{1,k}\right)=I_{k}^{-1}$, then the asymptotic variance of $N_{\tau_j,1}/\sqrt{\tau_j}$ for the Group ERADE with delayed and missing responses attains the Cramér-Rao lower bound
	$$\left(\left.\frac{\partial\rho}{\partial\mathbf{y}}\right|_{\boldsymbol{\Theta}}\right)' \mathrm{diag}\left(\left(\left(1-\beta_1\right)\nu I_1\right)^{-1},\left(\left(1-\beta_2\right)\left(1-\nu\right)I_2\right)^{-1}\right)\left.\frac{\partial\rho}{\partial\mathbf{y}}\right|_{\boldsymbol{\Theta}}$$
\end{theorem}

\begin{remark}
	When $\beta_k=0$, i.e., in the ideal scenario without missing responses, the matrix and asymptotic variance are given by
	$$\mathbf{V}=\mathrm{diag}\left(\frac{1}{\nu}{V}_1,\frac{1}{1-\nu}{V}_2\right)\quad \text{and} \quad \sigma^2 = \left(\left.\frac{\partial \rho}{\partial \mathbf{y}}\right|_{\boldsymbol{\Theta}}\right)'\mathbf{V} \left.\frac{\partial \rho}{\partial \mathbf{y}}\right|_{\boldsymbol{\Theta}}$$
	Correspondingly, under Conditions 1 and 2, and when Assumption 1 is satisfied, we have, as $\tau_j \to \infty$,
	$$\left|N_{\tau_j,1}-\tau_j\widehat{\rho}_j\right|=o_P\left(\sqrt{\tau_j}\right)\quad \text{and} \quad \left|N_{\tau_j,1}-\tau_j\widehat{\rho}_j\right|=O\left(\sqrt{\tau_j\log \log\tau_j}\right)$$
	Further, we have
	$$\sqrt{\tau_j}\left(N_{\tau_j,1}/\tau_j-\nu\right) \xrightarrow{D} \mathrm{N}\left(0,\sigma^{2}\right)$$
	and
	$$N_{\tau_j,1}/{\tau_j}-\nu=O\left(\sqrt{\tau_j^{-1} \log \log \tau_j}\right)$$
	The asymptotic variance of $N_{\tau_j,1}/\sqrt{\tau_j}$ for the Group ERADE without delayed and missing responses attains the Cramér-Rao lower bound
	$$\left(\left.\frac{\partial\rho}{\partial\mathbf{y}}\right|_{\boldsymbol{\Theta}}\right)' \mathrm{diag}((\nu I_1)^{-1},((1-\nu)I_2)^{-1})\left.\frac{\partial\rho}{\partial\mathbf{y}}\right|_{\boldsymbol{\Theta}}$$
\end{remark}

\subsection{Target Allocation Proportions}

In this section, we present several commonly used optimal allocation ratios, which can be broadly categorized into two main types: binary response and Gaussian response. We assume that the missing rates follow a Bernoulli distribution with mean $\beta_k$. To simplify the simulation, we also assume that the missing rate parameters are the same for all treatment groups, which means $\beta=\beta_1=\beta_2$.

In the case of binary response, we denote $P_1=1-Q_1$ and $P_2=1-Q_2$, where $P_1$ and $P_2$ represent the success probabilities of treatments $1$ and $2$, respectively. As defined earlier, natural choices for parameter estimation are
$$\widehat{P}_{j,k}=\frac{\sum_{m=1}^{j}\mathbf{X}_{k}^{\left(m\right)}\boldsymbol{\xi}_{k}^{\left(m\right)'}+0.5}{N_{\tau_j,k}+1}\quad\widehat{P}_{j,k}^{dm}=\dfrac{S_{\tau_j,k}+0.5}{T_{\tau_j,k}+1}\quad k=\left(1,2\right)$$
To target $\nu=\rho\left(P_1,P_2\right)$ as the allocation proportion, let $\widehat{\rho}_j=\rho\left(P_{j,1},P_{j,2}\right)$. According to Theorems 2.2 and 2.4, the variance of the allocation proportion achieves the Cramér-Rao lower bound under its respective allocations. We now list two speciﬁc allocation proportions:

(i) The RSIHR allocation proportion used by Rosenberger et al. \cite{rosenberger2001optimal},
\begin{align*}
	\nu_1=\frac{\sqrt{P_1}}{\sqrt{P_1}+\sqrt{P_2}} \quad \quad &\sigma_1^2=\frac{Q_2P_1^{1.5}+Q_1P_2^{1.5}}{4\sqrt{P_1P_2}\left(\sqrt{P_1}+\sqrt{P_2}\right)^3}\\
	&{\sigma_1^2}^{mis}=\frac{Q_2P_1^{1.5}+Q_1P_2^{1.5}}{4\left(1-\beta\right)\sqrt{P_1P_2}\left(\sqrt{P_1}+\sqrt{P_2}\right)^3}
\end{align*}

(ii) The Neyman allocation proportion discussed by Jennison and Turnbull \cite{jennison2000group},
\begin{align*}
	&\nu_2=\frac{\sqrt{P_1Q_1}}{\sqrt{P_1Q_1}+\sqrt{P_2Q_2}}\\ &\sigma_2^2=\frac{\left(P_1Q_1\right)^{1.5}\left(1-2P_2\right)^2+\left(P_2Q_2\right)^{1.5}\left(1-2P_1\right)^2}{4\sqrt{P_1Q_1P_2Q_2}\left(\sqrt{P_1Q_1}+\sqrt{P_2Q_2}\right)^3}\\
	&{\sigma_2^2}^{mis}=\frac{\left(P_1Q_1\right)^{1.5}\left(1-2P_2\right)^2+\left(P_2Q_2\right)^{1.5}\left(1-2P_1\right)^2}{4\left(1-\beta\right)\sqrt{P_1Q_1P_2Q_2}\left(\sqrt{P_1Q_1}+\sqrt{P_2Q_2}\right)^3}
\end{align*}

In the case of Gaussian response, We compare two treatments with responses $\xi_1\sim N\left(\mu_1,\pi_1^2\right)$ and $\xi_2\sim N\left(\mu_2,\pi_2^2\right)$, respectively. As defined earlier, natural choices for parameter estimation are
\begin{align*}
	&\widehat{\mu}_{j,k}=\frac{\sum_{m=1}^{j}\mathbf{X}_{k}^{\left(m\right)}\boldsymbol{\xi}_{k}^{\left(m\right)'}}{N_{\tau_j,k}}\\
	&\widehat{\pi}_{j,k}^2=\frac{\sum_{m=1}^{j}\sum_{i=1}^{\chi_m}X_{i,k}^{\left(m\right)}\left({\xi}_{i,k}^{\left(m\right)}-\widehat{\mu}_{j,k}\right)^2}{N_{\tau_j,k}}\\
	&\widehat{\mu}_{j,k}^{dm}=\frac{S_{\tau_j,k}}{T_{\tau_j,k}}\\
	&\widehat{\pi}_{j,k}^{2dm}=\frac{\sum_{m=1}^{j}\sum_{l=0}^{j-m}\sum_{i=1}^{\chi_m}M_{i,k}(m,l)\big(1-I_{mis,k}^{(m)}(i)\big)X_{i,k}^{(m)}\big({\xi}_{i,k}^{\left(m\right)}-\widehat{\mu}_{j,k}\big)^2}{T_{\tau_j,k}}
\end{align*}

We use the estimator $\widehat{\rho}_j=\rho\left(\widehat{\mu}_{j,1},\widehat{\mu}_{j,2},\widehat{\pi}_{j,1},\widehat{\pi}_{j,2}\right)$ to target a given allocation proportion $\nu$. According to Theorems 2.2 and 2.4, the variance of the allocation proportion achieves the Cramér-Rao lower bound under its respective allocations. We now list two speciﬁc allocation proportions:

(i) The optimal allocation proportion discussed by Zhang and Rosenberger \cite{zhang2006response},
\begin{align*}
	\nu_3=\frac{\pi_1\sqrt{\mu_2}}{\pi_1\sqrt{\mu_2}+\pi_2\sqrt{\mu_1}}\quad \quad &\sigma_3^2=\frac{\pi_1\pi_2\sqrt{\mu_1}\sqrt{\mu_2}}{2\left(\pi_1\sqrt{\mu_2}+\pi_2\sqrt{\mu_1}\right)^2}\\
	&{\sigma_3^2}^{mis}=\frac{\pi_1\pi_2\sqrt{\mu_1}\sqrt{\mu_2}}{2\left(1-\beta\right)\left(\pi_1\sqrt{\mu_2}+\pi_2\sqrt{\mu_1}\right)^2}
\end{align*}

(ii) The Neyman allocation proportion given by Jennison and Turnbull \cite{jennison2000group},
\begin{align*}
	\nu_4=\frac{\pi_1}{\pi_1+\pi_2}\quad \quad &\sigma_4^2=\frac{\pi_1\pi_2}{2\left(\pi_1+\pi_2\right)^2}\\
	&{\sigma_4^2}^{mis}=\frac{\pi_1\pi_2}{2\left(1-\beta\right)\left(\pi_1+\pi_2\right)^2}
\end{align*}

\begin{remark}
	All of them are optimal allocation proportions. The RSIHR aims to minimize the expected number of treatment failures in the trial while maintaining a fixed variance of the test statistic, whereas the Neyman allocation proportion maximizes the power when the sample size is fixed.
\end{remark}

\section{Simulation Study}\label{sec3}

In this section, we conducted a series of numerical simulation studies and compared it with the ERADE proposed by Hu, Zhang, and He \cite{hu2009efficient} to evaluate the empirical performance of Group ERADE proposed in this paper and its improved variant under delayed and missing responses. The primary aim is to verify that Group ERADE not only aligns more closely with the actual operational workflow of clinical trials, but also retains the advantages of ERADE. We selected the mean and standard deviation of the realized allocation proportion and the cumulative number of treatment failures in the trial as the core comparison metrics. 

In the simulation setup, restricted randomization was first applied to allocate 10 patients to each treatment arm as the initial stage. Subsequently, the remaining 180 patients were sequentially allocated according to the design. The group size follows a Poisson distribution with mean $\lambda=10$, and the response delay time for each treatment is assumed to follow an exponential distribution with mean $1/\lambda_k\ \left(k=1,2\right)$. In this study, we set $\lambda_k=10$.

\subsection{The Case of Binary Response}

In the binary response setting, we follow the parameter configuration in Hu, Zhang, and He \cite{hu2009efficient}. To control the length of the paper, this article only reports the results obtained using the RSIHR target allocation proportion (see Tables 1 and 2, and the results under the Neyman allocation proportion are detailed in the appendix. All simulations are based on 5000 independent replications, with a total sample size of $\tau_j=200$.

The results in Tables 1 and 2 show that under various combinations of $P_1$ and $P_2$, the Group ERADE, including its improved variant under delayed and missing responses, can approximate the target allocation proportion well when $\alpha=2/3$. The variance of the simulated allocation proportion is generally consistent with the theoretical asymptotic variance, and the cumulative number of treatment failures remains at a low level. Compared with Group DBCD method with $\gamma=2$, Group ERADE exhibits a lower simulation variance. In addition, we also find that under the same missing rate condition, when the success probabilities $P_1$ and $P_2$ of the two treatments are both relatively small, the performance of Group ERADE is relatively better.

\begin{table}[h]
	\centering
	\small
	\caption{Results for the binary case with RSIHR target and no delayed and missing responses.}
	\resizebox{\textwidth}{!}{
		\begin{tabular}{cccccccc}
			\hline
			& & \multicolumn{2}{c}{ERADE} & \multicolumn{2}{c}{Group ERADE} & \multicolumn{2}{c}{Group DBCD} \\
			\cline{3-8} 
			$P_1$, $P_2$​ & \makecell[c]{$\nu$\\$\sigma/\sqrt{\tau_j}$} & \makecell[c]{$\mathbb{E}\left({N_{\tau_j,1}}/{\tau_j}\right)$\\(SD)} & \makecell[c]{$\mathbb{E}\left(\text{No.Fail}\right)$\\(SD)} & \makecell[c]{$\mathbb{E}\left({N_{\tau_j,1}}/{\tau_j}\right)$\\(SD)} & \makecell[c]{$\mathbb{E}\left(\text{No.Fail}\right)$\\(SD)} & \makecell[c]{$\mathbb{E}\left({N_{\tau_j,1}}/{\tau_j}\right)$\\(SD)} & \makecell[c]{$\mathbb{E}\left(\text{No.Fail}\right)$\\(SD)}\\\hline
			0.9, 0.7 & \makecell[c]{0.531\\(0.009)} & \makecell[c]{0.531\\(0.011)} & \makecell[c]{0.194\\(0.026)} & \makecell[c]{0.530\\(0.013)} & \makecell[c]{0.194\\(0.027)} & \makecell[c]{0.531\\(0.019)} & \makecell[c]{0.194\\(0.026)} \\\hline
			0.9, 0.5 & \makecell[c]{0.573\\(0.014)} & \makecell[c]{0.572\\(0.015)} & \makecell[c]{0.271\\(0.025)} & \makecell[c]{0.571\\(0.017)} & \makecell[c]{0.271\\(0.026)} & \makecell[c]{0.573\\(0.023)} & \makecell[c]{0.271\\(0.026)} \\\hline
			0.8, 0.8 & \makecell[c]{0.500\\(0.009)} & \makecell[c]{0.500\\(0.011)} & \makecell[c]{0.200\\(0.028)} & \makecell[c]{0.500\\(0.013)} & \makecell[c]{0.200\\(0.029)} & \makecell[c]{0.499\\(0.019)} & \makecell[c]{0.200\\(0.028)} \\\hline
			0.8, 0.6 & \makecell[c]{0.536\\(0.012)} & \makecell[c]{0.535\\(0.013)} & \makecell[c]{0.293\\(0.031)} & \makecell[c]{0.535\\(0.016)} & \makecell[c]{0.293\\(0.031)} & \makecell[c]{0.536\\(0.021)} & \makecell[c]{0.292\\(0.031)} \\\hline
			0.7, 0.5 & \makecell[c]{0.542\\(0.015)} & \makecell[c]{0.542\\(0.017)} & \makecell[c]{0.392\\(0.033)} & \makecell[c]{0.541\\(0.018)} & \makecell[c]{0.392\\(0.034)} & \makecell[c]{0.542\\(0.024)} & \makecell[c]{0.391\\(0.033)} \\\hline
			0.6, 0.4 & \makecell[c]{0.551\\(0.019)} & \makecell[c]{0.550\\(0.020)} & \makecell[c]{0.490\\(0.034)} & \makecell[c]{0.549\\(0.022)} & \makecell[c]{0.491\\(0.034)} & \makecell[c]{0.551\\(0.027)} & \makecell[c]{0.490\\(0.034)} \\\hline
			0.5, 0.5 & \makecell[c]{0.500\\(0.018)} & \makecell[c]{0.500\\(0.019)} & \makecell[c]{0.500\\(0.036)} & \makecell[c]{0.500\\(0.020)} & \makecell[c]{0.501\\(0.035)} & \makecell[c]{0.500\\(0.026)} & \makecell[c]{0.500\\(0.035)} \\\hline
			0.4, 0.3 & \makecell[c]{0.536\\(0.025)} & \makecell[c]{0.535\\(0.027)} & \makecell[c]{0.647\\(0.033)} & \makecell[c]{0.536\\(0.028)} & \makecell[c]{0.646\\(0.034)} & \makecell[c]{0.538\\(0.034)} & \makecell[c]{0.646\\(0.034)} \\\hline
			0.2, 0.2 & \makecell[c]{0.500\\(0.035)} & \makecell[c]{0.500\\(0.040)} & \makecell[c]{0.800\\(0.028)} & \makecell[c]{0.500\\(0.041)} & \makecell[c]{0.800\\(0.029)} & \makecell[c]{0.499\\(0.045)} & \makecell[c]{0.800\\(0.028)} \\\hline
		\end{tabular}
	}
\end{table}

\begin{table}[h]
	\centering
	\small
	\caption{Delayed and missing response: results for the binary case with RSIHR target.}
	\resizebox{\textwidth}{!}{
		\begin{tabular}{cccccccc}
			\hline
			& & \multicolumn{2}{c}{ERADE} & \multicolumn{2}{c}{Group ERADE} & \multicolumn{2}{c}{Group DBCD} \\
			\cline{3-8} 
			$\beta$, $P_1$, $P_2$​ & \makecell[c]{$\nu$\\$\sigma/\sqrt{\tau_j}$} & \makecell[c]{$\mathbb{E}\left({N_{\tau_j,1}}/{\tau_j}\right)$\\(SD)} & \makecell[c]{$\mathbb{E}\left(\text{No.Fail}\right)$\\(SD)} & \makecell[c]{$\mathbb{E}\left({N_{\tau_j,1}}/{\tau_j}\right)$\\(SD)} & \makecell[c]{$\mathbb{E}\left(\text{No.Fail}\right)$\\(SD)} & \makecell[c]{$\mathbb{E}\left({N_{\tau_j,1}}/{\tau_j}\right)$\\(SD)} & \makecell[c]{$\mathbb{E}\left(\text{No.Fail}\right)$\\(SD)}\\\hline
			0.0 0.9, 0.7 & \makecell[c]{0.531\\(0.009)} & \makecell[c]{0.530\\(0.012)} & \makecell[c]{0.194\\(0.026)} & \makecell[c]{0.531\\(0.014)} & \makecell[c]{0.193\\(0.026)} & \makecell[c]{0.531\\(0.020)} & \makecell[c]{0.193\\(0.026)} \\\hline
			0.0 0.9, 0.5 & \makecell[c]{0.573\\(0.014)} & \makecell[c]{0.571\\(0.016)} & \makecell[c]{0.270\\(0.025)} & \makecell[c]{0.571\\(0.017)} & \makecell[c]{0.271\\(0.025)} & \makecell[c]{0.573\\(0.023)} & \makecell[c]{0.270\\(0.026)} \\\hline
			0.0 0.8, 0.8 & \makecell[c]{0.500\\(0.009)} & \makecell[c]{0.500\\(0.011)} & \makecell[c]{0.200\\(0.029)} & \makecell[c]{0.500\\(0.013)} & \makecell[c]{0.200\\(0.028)} & \makecell[c]{0.500\\(0.019)} & \makecell[c]{0.200\\(0.028)} \\\hline
			0.0 0.8, 0.6 & \makecell[c]{0.536\\(0.012)} & \makecell[c]{0.535\\(0.015)} & \makecell[c]{0.293\\(0.031)} & \makecell[c]{0.535\\(0.016)} & \makecell[c]{0.293\\(0.030)} & \makecell[c]{0.536\\(0.022)} & \makecell[c]{0.293\\(0.031)} \\\hline
			0.0 0.7, 0.5 & \makecell[c]{0.542\\(0.015)} & \makecell[c]{0.540\\(0.018)} & \makecell[c]{0.392\\(0.033)} & \makecell[c]{0.541\\(0.018)} & \makecell[c]{0.391\\(0.033)} & \makecell[c]{0.543\\(0.024)} & \makecell[c]{0.392\\(0.034)} \\\hline
			0.0 0.6, 0.4 & \makecell[c]{0.551\\(0.019)} & \makecell[c]{0.549\\(0.022)} & \makecell[c]{0.490\\(0.034)} & \makecell[c]{0.549\\(0.022)} & \makecell[c]{0.490\\(0.034)} & \makecell[c]{0.551\\(0.028)} & \makecell[c]{0.490\\(0.034)} \\\hline
			0.0 0.5, 0.5 & \makecell[c]{0.500\\(0.018)} & \makecell[c]{0.499\\(0.021)} & \makecell[c]{0.500\\(0.035)} & \makecell[c]{0.500\\(0.021)} & \makecell[c]{0.501\\(0.035)} & \makecell[c]{0.500\\(0.026)} & \makecell[c]{0.500\\(0.035)} \\\hline
			0.0 0.4, 0.3 & \makecell[c]{0.536\\(0.025)} & \makecell[c]{0.534\\(0.031)} & \makecell[c]{0.647\\(0.034)} & \makecell[c]{0.536\\(0.028)} & \makecell[c]{0.646\\(0.033)} & \makecell[c]{0.537\\(0.033)} & \makecell[c]{0.647\\(0.034)} \\\hline
			0.0 0.2, 0.2 & \makecell[c]{0.500\\(0.035)} & \makecell[c]{0.497\\(0.045)} & \makecell[c]{0.800\\(0.029)} & \makecell[c]{0.500\\(0.043)} & \makecell[c]{0.800\\(0.028)} & \makecell[c]{0.500\\(0.047)} & \makecell[c]{0.800\\(0.029)} \\\hline
			0.1 0.9, 0.7 & \makecell[c]{0.531\\(0.010)} & \makecell[c]{0.530\\(0.017)} & \makecell[c]{0.193\\(0.028)} & \makecell[c]{0.531\\(0.019)} & \makecell[c]{0.194\\(0.028)} & \makecell[c]{0.531\\(0.023)} & \makecell[c]{0.194\\(0.028)} \\\hline
			0.1 0.9, 0.5 & \makecell[c]{0.573\\(0.015)} & \makecell[c]{0.571\\(0.021)} & \makecell[c]{0.271\\(0.027)} & \makecell[c]{0.571\\(0.021)} & \makecell[c]{0.272\\(0.027)} & \makecell[c]{0.573\\(0.027)} & \makecell[c]{0.271\\(0.028)} \\\hline
			0.1 0.8, 0.8 & \makecell[c]{0.500\\(0.009)} & \makecell[c]{0.500\\(0.017)} & \makecell[c]{0.200\\(0.030)} & \makecell[c]{0.500\\(0.018)} & \makecell[c]{0.200\\(0.030)} & \makecell[c]{0.500\\(0.023)} & \makecell[c]{0.200\\(0.030)} \\\hline
			0.1 0.8, 0.6 & \makecell[c]{0.536\\(0.013)} & \makecell[c]{0.535\\(0.019)} & \makecell[c]{0.293\\(0.033)} & \makecell[c]{0.535\\(0.020)} & \makecell[c]{0.293\\(0.032)} & \makecell[c]{0.536\\(0.025)} & \makecell[c]{0.293\\(0.033)} \\\hline
			0.1 0.7, 0.5 & \makecell[c]{0.542\\(0.016)} & \makecell[c]{0.540\\(0.022)} & \makecell[c]{0.392\\(0.035)} & \makecell[c]{0.541\\(0.023)} & \makecell[c]{0.392\\(0.035)} & \makecell[c]{0.542\\(0.027)} & \makecell[c]{0.391\\(0.036)} \\\hline
			0.1 0.6, 0.4 & \makecell[c]{0.551\\(0.020)} & \makecell[c]{0.549\\(0.026)} & \makecell[c]{0.490\\(0.036)} & \makecell[c]{0.550\\(0.026)} & \makecell[c]{0.490\\(0.036)} & \makecell[c]{0.552\\(0.031)} & \makecell[c]{0.490\\(0.036)} \\\hline
			0.1 0.5, 0.5 & \makecell[c]{0.500\\(0.019)} & \makecell[c]{0.499\\(0.025)} & \makecell[c]{0.500\\(0.037)} & \makecell[c]{0.500\\(0.024)} & \makecell[c]{0.500\\(0.037)} & \makecell[c]{0.501\\(0.030)} & \makecell[c]{0.500\\(0.037)} \\\hline
			0.1 0.4, 0.3 & \makecell[c]{0.536\\(0.026)} & \makecell[c]{0.535\\(0.034)} & \makecell[c]{0.646\\(0.035)} & \makecell[c]{0.536\\(0.032)} & \makecell[c]{0.647\\(0.035)} & \makecell[c]{0.537\\(0.038)} & \makecell[c]{0.646\\(0.035)} \\\hline
			0.1 0.2, 0.2 & \makecell[c]{0.500\\(0.037)} & \makecell[c]{0.497\\(0.052)} & \makecell[c]{0.800\\(0.029)} & \makecell[c]{0.500\\(0.046)} & \makecell[c]{0.800\\(0.030)} & \makecell[c]{0.501\\(0.052)} & \makecell[c]{0.800\\(0.030)} \\\hline
			0.2 0.9, 0.7 & \makecell[c]{0.531\\(0.010)} & \makecell[c]{0.530\\(0.022)} & \makecell[c]{0.195\\(0.030)} & \makecell[c]{0.531\\(0.023)} & \makecell[c]{0.194\\(0.029)} & \makecell[c]{0.532\\(0.027)} & \makecell[c]{0.194\\(0.029)} \\\hline
			0.2 0.9, 0.5 & \makecell[c]{0.573\\(0.015)} & \makecell[c]{0.571\\(0.025)} & \makecell[c]{0.271\\(0.030)} & \makecell[c]{0.571\\(0.026)} & \makecell[c]{0.271\\(0.030)} & \makecell[c]{0.574\\(0.030)} & \makecell[c]{0.270\\(0.030)} \\\hline
			0.2 0.8, 0.8 & \makecell[c]{0.500\\(0.010)} & \makecell[c]{0.499\\(0.021)} & \makecell[c]{0.200\\(0.032)} & \makecell[c]{0.500\\(0.023)} & \makecell[c]{0.200\\(0.032)} & \makecell[c]{0.500\\(0.027)} & \makecell[c]{0.200\\(0.032)} \\\hline
			0.2 0.8, 0.6 & \makecell[c]{0.536\\(0.014)} & \makecell[c]{0.535\\(0.024)} & \makecell[c]{0.293\\(0.035)} & \makecell[c]{0.535\\(0.025)} & \makecell[c]{0.293\\(0.034)} & \makecell[c]{0.536\\(0.028)} & \makecell[c]{0.293\\(0.035)} \\\hline
			0.2 0.7, 0.5 & \makecell[c]{0.542\\(0.017)} & \makecell[c]{0.540\\(0.026)} & \makecell[c]{0.392\\(0.037)} & \makecell[c]{0.540\\(0.027)} & \makecell[c]{0.392\\(0.037)} & \makecell[c]{0.542\\(0.032)} & \makecell[c]{0.391\\(0.038)} \\\hline
			0.2 0.6, 0.4 & \makecell[c]{0.551\\(0.021)} & \makecell[c]{0.549\\(0.031)} & \makecell[c]{0.490\\(0.038)} & \makecell[c]{0.549\\(0.030)} & \makecell[c]{0.490\\(0.038)} & \makecell[c]{0.552\\(0.035)} & \makecell[c]{0.490\\(0.038)} \\\hline
			0.2 0.5, 0.5 & \makecell[c]{0.500\\(0.020)} & \makecell[c]{0.499\\(0.030)} & \makecell[c]{0.499\\(0.039)} & \makecell[c]{0.501\\(0.029)} & \makecell[c]{0.500\\(0.040)} & \makecell[c]{0.501\\(0.034)} & \makecell[c]{0.500\\(0.039)} \\\hline
			0.2 0.4, 0.3 & \makecell[c]{0.536\\(0.028)} & \makecell[c]{0.535\\(0.037)} & \makecell[c]{0.647\\(0.038)} & \makecell[c]{0.536\\(0.037)} & \makecell[c]{0.647\\(0.038)} & \makecell[c]{0.537\\(0.042)} & \makecell[c]{0.646\\(0.037)} \\\hline
			0.2 0.2, 0.2 & \makecell[c]{0.500\\(0.040)} & \makecell[c]{0.496\\(0.057)} & \makecell[c]{0.800\\(0.031)} & \makecell[c]{0.499\\(0.052)} & \makecell[c]{0.800\\(0.032)} & \makecell[c]{0.502\\(0.056)} & \makecell[c]{0.800\\(0.031)} \\\hline
		\end{tabular}
	}
\end{table}

\subsection{The Case of Continuous Response}

Next, we further evaluate the performance of Group ERADE for continuous response variables. Taking the optimal allocation proportion proposed by Zhang and Rosenberger  \cite{zhang2006response} as the target, the relevant numerical results are summarized in Tables 3 and 4. We examined multiple sets of different parameter combinations $(\mu_1,\pi_1,\mu_2,\pi_2)$ to comprehensively reflect the design performance under different treatment effects and degrees of variation.

\begin{table}[h]
	\centering
	\small
	\caption{Results for the continuous case with Zhang and Rosenberger's target.}
	\resizebox{\textwidth}{!}{
		\begin{tabular}{cccccccc}
			\hline
			& & \multicolumn{2}{c}{ERADE} & \multicolumn{2}{c}{Group ERADE} & \multicolumn{2}{c}{Group DBCD} \\
			\cline{3-8} 
			$\mu_1$,$\pi_1$,$\mu_2$,$\pi_2$​ & \makecell[c]{$\nu$\\$\sigma/\sqrt{\tau_j}$} & \makecell[c]{$\mathbb{E}\left({N_{\tau_j,1}}/{\tau_j}\right)$\\(SD)} & \makecell[c]{$\mathbb{E}\left(\text{Obs.Rsp}\right)$\\(SD)} & \makecell[c]{$\mathbb{E}\left({N_{\tau_j,1}}/{\tau_j}\right)$\\(SD)} & \makecell[c]{$\mathbb{E}\left(\text{Obs.Rsp}\right)$\\(SD)} & \makecell[c]{$\mathbb{E}\left({N_{\tau_j,1}}/{\tau_j}\right)$\\(SD)} & \makecell[c]{$\mathbb{E}\left(\text{Obs.Rsp}\right)$\\(SD)}\\\hline
			13,4,15,2.5 & \makecell[c]{0.632\\(0.024)} & \makecell[c]{0.631\\(0.026)} & \makecell[c]{13.732\\(0.257)} & \makecell[c]{0.629\\(0.027)} & \makecell[c]{13.740\\(0.261)} & \makecell[c]{0.634\\(0.032)} & \makecell[c]{13.732\\(0.254)} \\\hline
			13,2.5,15,4 & \makecell[c]{0.402\\(0.025)} & \makecell[c]{0.403\\(0.026)} & \makecell[c]{14.196\\(0.251)} & \makecell[c]{0.404\\(0.027)} & \makecell[c]{14.188\\(0.250)} & \makecell[c]{0.401\\(0.032)} & \makecell[c]{14.203\\(0.252)} \\\hline
			15,4,17,2.5 & \makecell[c]{0.630\\(0.024)} & \makecell[c]{0.629\\(0.026)} & \makecell[c]{15.737\\(0.256)} & \makecell[c]{0.627\\(0.027)} & \makecell[c]{15.742\\(0.261)} & \makecell[c]{0.632\\(0.032)} & \makecell[c]{15.727\\(0.263)} \\\hline
			15,2.5,17,4 & \makecell[c]{0.400\\(0.024)} & \makecell[c]{0.400\\(0.027)} & \makecell[c]{16.198\\(0.249)} & \makecell[c]{0.401\\(0.028)} & \makecell[c]{16.203\\(0.251)} & \makecell[c]{0.397\\(0.033)} & \makecell[c]{16.204\\(0.253)} \\\hline
		\end{tabular}
	}
\end{table}

\begin{table}[h]
	\centering
	\small
	\caption{Delayed and missing response: results for the continuous case with Zhang and Rosenberger's target.}
	\resizebox{\textwidth}{!}{
		\begin{tabular}{cccccccc}
			\hline
			& & \multicolumn{2}{c}{ERADE} & \multicolumn{2}{c}{Group ERADE} & \multicolumn{2}{c}{Group DBCD} \\
			\cline{3-8} 
			$\beta$,$\mu_1$,$\pi_1$,$\mu_2$,$\pi_2$​ & \makecell[c]{$\nu$\\$\sigma/\sqrt{\tau_j}$} & \makecell[c]{$\mathbb{E}\left({N_{\tau_j,1}}/{\tau_j}\right)$\\(SD)} & \makecell[c]{$\mathbb{E}\left(\text{Obs.Rsp}\right)$\\(SD)} & \makecell[c]{$\mathbb{E}\left({N_{\tau_j,1}}/{\tau_j}\right)$\\(SD)} & \makecell[c]{$\mathbb{E}\left(\text{Obs.Rsp}\right)$\\(SD)} & \makecell[c]{$\mathbb{E}\left({N_{\tau_j,1}}/{\tau_j}\right)$\\(SD)} & \makecell[c]{$\mathbb{E}\left(\text{Obs.Rsp}\right)$\\(SD)}\\\hline
			0.0,13,4,15,2.5 & \makecell[c]{0.632\\(0.024)} & \makecell[c]{0.630\\(0.026)} & \makecell[c]{13.738\\(0.258)} & \makecell[c]{0.629\\(0.028)} & \makecell[c]{13.742\\(0.261)} & \makecell[c]{0.633\\(0.032)} & \makecell[c]{13.736\\(0.258)} \\\hline
			0.0,13,2.5,15,4 & \makecell[c]{0.402\\(0.025)} & \makecell[c]{0.403\\(0.027)} & \makecell[c]{14.191\\(0.246)} & \makecell[c]{0.404\\(0.029)} & \makecell[c]{14.191\\(0.253)} & \makecell[c]{0.401\\(0.034)} & \makecell[c]{14.200\\(0.253)} \\\hline
			0.0,15,4,17,2.5 & \makecell[c]{0.630\\(0.024)} & \makecell[c]{0.629\\(0.026)} & \makecell[c]{15.742\\(0.258)} & \makecell[c]{0.627\\(0.027)} & \makecell[c]{15.745\\(0.256)} & \makecell[c]{0.632\\(0.033)} & \makecell[c]{15.743\\(0.262)} \\\hline
			0.0,15,2.5,17,4 & \makecell[c]{0.400\\(0.024)} & \makecell[c]{0.401\\(0.026)} & \makecell[c]{16.201\\(0.245)} & \makecell[c]{0.402\\(0.028)} & \makecell[c]{16.195\\(0.252)} & \makecell[c]{0.398\\(0.033)} & \makecell[c]{16.206\\(0.250)} \\\hline
			0.2,13,4,15,2.5 & \makecell[c]{0.632\\(0.027)} & \makecell[c]{0.630\\(0.034)} & \makecell[c]{13.737\\(0.292)} & \makecell[c]{0.629\\(0.036)} & \makecell[c]{13.740\\(0.293)} & \makecell[c]{0.634\\(0.040)} & \makecell[c]{13.734\\(0.291)} \\\hline
			0.2,13,2.5,15,4 & \makecell[c]{0.402\\(0.027)} & \makecell[c]{0.402\\(0.035)} & \makecell[c]{14.191\\(0.278)} & \makecell[c]{0.404\\(0.037)} & \makecell[c]{14.193\\(0.282)} & \makecell[c]{0.400\\(0.041)} & \makecell[c]{14.198\\(0.283)} \\\hline
			0.2,15,4,17,2.5 & \makecell[c]{0.630\\(0.027)} & \makecell[c]{0.629\\(0.034)} & \makecell[c]{15.743\\(0.288)} & \makecell[c]{0.627\\(0.036)} & \makecell[c]{15.740\\(0.293)} & \makecell[c]{0.632\\(0.041)} & \makecell[c]{15.734\\(0.293)} \\\hline
			0.2,15,2.5,17,4 & \makecell[c]{0.400\\(0.027)} & \makecell[c]{0.400\\(0.034)} & \makecell[c]{16.202\\(0.278)} & \makecell[c]{0.402\\(0.036)} & \makecell[c]{16.192\\(0.280)} & \makecell[c]{0.398\\(0.041)} & \makecell[c]{16.201\\(0.279)} \\\hline
			0.3,13,4,15,2.5 & \makecell[c]{0.632\\(0.029)} & \makecell[c]{0.631\\(0.040)} & \makecell[c]{13.735\\(0.315)} & \makecell[c]{0.629\\(0.041)} & \makecell[c]{13.743\\(0.312)} & \makecell[c]{0.634\\(0.046)} & \makecell[c]{13.733\\(0.317)} \\\hline
			0.3,13,2.5,15,4 & \makecell[c]{0.402\\(0.029)} & \makecell[c]{0.401\\(0.041)} & \makecell[c]{14.192\\(0.306)} & \makecell[c]{0.404\\(0.041)} & \makecell[c]{14.196\\(0.297)} & \makecell[c]{0.400\\(0.046)} & \makecell[c]{14.198\\(0.304)} \\\hline
			0.3,15,4,17,2.5 & \makecell[c]{0.630\\(0.029)} & \makecell[c]{0.629\\(0.040)} & \makecell[c]{15.739\\(0.309)} & \makecell[c]{0.628\\(0.041)} & \makecell[c]{15.741\\(0.315)} & \makecell[c]{0.632\\(0.046)} & \makecell[c]{15.733\\(0.317)} \\\hline
			0.3,15,2.5,17,4 & \makecell[c]{0.400\\(0.029)} & \makecell[c]{0.401\\(0.040)} & \makecell[c]{16.200\\(0.303)} & \makecell[c]{0.402\\(0.041)} & \makecell[c]{16.190\\(0.304)} & \makecell[c]{0.398\\(0.047)} & \makecell[c]{16.202\\(0.304)} \\\hline
		\end{tabular}
	}
\end{table}

Under the continuous response setting, we use the sample mean of the observed responses as the primary efficacy endpoint, rather than the failure rate used in the binary case. The simulation results show that there is only a small deviation between the actual allocation proportion and the target value for both Group ERADE and the improved version designed for delayed and missing responses. In addition, Group ERADE and ERADE show a high degree of similarity in terms of the mean and standard deviation of the allocation proportion, indicating that they have comparable statistical precision. These results fully verify the robustness and reliability of the proposed method in actual clinical trials, and the method is especially suitable for scenarios with varying missing rates, delay times, and group sizes. Compared with Group DBCD, Group ERADE exhibits lower simulation variance. The simulation results using the Neyman allocation proportion as the target are detailed in the appendix, where similar results are observed.

\subsection{A Redesigned Clinical Trial}

To verify the practical application value of the proposed method, we referred to a clinical trial conducted by Dworkin et al. \cite{dworkin2003pregabalin} evaluating pregabalin for the treatment of post-herpetic neuralgia (PHN). This study was a multi-center, parallel-group, double-blind, randomized controlled trial. A total of $173$ PHN patients were enrolled over an $8$-week trial period. Among them, $89$ patients were assigned to the pregabalin group and $84$ to the placebo group. The results showed that $63\%$ of patients in the pregabalin group reported a pain reduction of $\ge30\%$, while the proportion in the placebo group was $25\%$; for the endpoint of pain reduction $\ge50\%$, the proportions in the two groups were $50\%$ and $20\%$ respectively.

Based on the above real research data, we carried out a retrospective re-design. We adopted the Group ERADE proposed in this paper (with $\alpha=2/3$) and its improved version under delayed and missing responses. In the re-designed simulated trial, we used the real remission rates reported in the original literature as the basis for data generation, and assumed that these real parameters were unknown throughout the trial process and were only used to simulate patient-level response results. The enrollment process was consistent with the original trial: $173$ PHN patients (defined as patients with pain lasting more than $3$ months after the healing of herpes zoster skin lesions) were recruited over $8$ weeks. Referring to the original literature, the missing rates for both treatment types were set as $\beta_1=\beta_2=0.24$. Two enrollment rhythms were also considered: (i) updating the grouping every two days, corresponding to $\lambda = 6.18$; (ii) updating the grouping every four days, corresponding to $\lambda = 12.36$. The target allocation proportion was determined using the RSIHR allocation proportion.

In the initial stage, the first $10$ patients were assigned to different treatments using a simple randomization method to obtain initial parameter estimates. Subsequently, multiple adaptive designs were successively applied to allocate the remaining $163$ patients, including Group ERADE proposed in this paper, DL, and Group DBCD. The comparison results of the above different design schemes are summarized in Tables 5 and 6.

\begin{table}[h]
	\centering
	\small
	\caption{The redesigned trial.(reduction in pain of $30\%$ or more)}
	\resizebox{\textwidth}{!}{
		\begin{tabular}{cccccccc}
			\hline
			Method & Update cycle & $\nu$ & Effective $n$ & Patients in preganalin & \makecell[c]{\text{Number of effective }\\\text{pain relief cases}} \\\hline
			Original & —— & —— & 173 & 89 & 77 \\\hline
			DL & —— & 0.670 & 173 & 114.3(5.3) & 86.7(7.2) \\\hline
			DBCD & —— & 0.614 & 173 & 105.2(6.3) & 82.3(6.1) \\\hline
			ERADE & —— & 0.614 & 173 & 106.2(5.3) & 83.6(6.0) \\\hline
			ERADE(DM) & —— & 0.614 & 131 & 80.7(6.7) & 63.3(5.4) \\\hline
			\makecell[c]{\text{Group}\\\text{DBCD}} & \makecell[c]{\text{2 days}\\\text{4 days}} & \makecell[c]{0.614\\0.614} & \makecell[c]{173\\173} & \makecell[c]{106.5(6.4)\\106.5(6.4)} & \makecell[c]{83.9(6.0)\\83.9(6.0)} \\\hline
			\makecell[c]{\text{Group}\\\text{DBCD(DM)}} & \makecell[c]{\text{2 days}\\\text{4 days}} & \makecell[c]{0.614\\0.614} & \makecell[c]{131\\131} & \makecell[c]{81.1(7.0)\\81.1(7.2)} & \makecell[c]{63.3(5.4)\\63.5(5.4)} \\\hline
			\makecell[c]{\text{Group}\\\text{ERADE}} & \makecell[c]{\text{2 days}\\\text{4 days}} & \makecell[c]{0.614\\0.614} & \makecell[c]{173\\173} & \makecell[c]{105.9(5.4)\\105.6(5.5)} & \makecell[c]{83.6(5.9)\\83.4(6.0)} \\\hline
			\makecell[c]{\text{Group}\\\text{ERADE(DM)}} & \makecell[c]{\text{2 days}\\\text{4 days}} & \makecell[c]{0.614\\0.614} & \makecell[c]{131\\131} & \makecell[c]{80.5(6.6)\\80.4(6.6)} & \makecell[c]{63.3(5.3)\\63.3(5.3)} \\\hline
		\end{tabular}
	}
\end{table}

\begin{table}[h]
	\centering
	\small
	\caption{The redesigned trial.(reduction in pain of $50\%$ or more)}
	\resizebox{\textwidth}{!}{
		\begin{tabular}{cccccccc}
			\hline
			Method & Update cycle & $\nu$ & Effective $n$ & Patients in preganalin & \makecell[c]{\text{Number of effective }\\\text{pain relief cases}} \\\hline
			Original & —— & —— & 173 & 89 & 61 \\\hline
			DL & —— & 0.615 & 173 & 105.8(4.7) & 66.4(6.8) \\\hline
			DBCD & —— & 0.613 & 173 & 105.0(7.3) & 65.5(6.2) \\\hline
			ERADE & —— & 0.614 & 173 & 106.0(6.2) & 66.5(6.1) \\\hline
			ERADE(DM) & —— & 0.614 & 131 & 80.4(7.6) & 50.4(5.4) \\\hline
			\makecell[c]{\text{Group}\\\text{DBCD}} & \makecell[c]{\text{2 days}\\\text{4 days}} & \makecell[c]{0.614\\0.614} & \makecell[c]{173\\173} & \makecell[c]{106.6(7.5)\\106.7(7.3)} & \makecell[c]{66.6(6.2)\\66.7(6.2)} \\\hline
			\makecell[c]{\text{Group}\\\text{DBCD(DM)}} & \makecell[c]{\text{2 days}\\\text{4 days}} & \makecell[c]{0.614\\0.614} & \makecell[c]{131\\131} & \makecell[c]{81.1(7.8)\\81.1(8.0)} & \makecell[c]{50.5(5.4)\\50.4(5.5)} \\\hline
			\makecell[c]{\text{Group}\\\text{ERADE}} & \makecell[c]{\text{2 days}\\\text{4 days}} & \makecell[c]{0.614\\0.614} & \makecell[c]{173\\173} & \makecell[c]{105.9(6.4)\\105.5(6.4)} & \makecell[c]{66.5(6.1)\\66.3(6.0)} \\\hline
			\makecell[c]{\text{Group}\\\text{ERADE(DM)}} & \makecell[c]{\text{2 days}\\\text{4 days}} & \makecell[c]{0.614\\0.614} & \makecell[c]{131\\131} & \makecell[c]{80.8(7.3)\\80.3(7.4)} & \makecell[c]{50.5(5.4)\\50.2(5.4)} \\\hline
		\end{tabular}
	}
\end{table}

Based on the results of the above redesigned experiments, if the same sample size as the original study is maintained, all the methods proposed in this paper can obtain more cases of effective pain relief with the same patient scale. Specifically, Group ERADE proposed by us and the original ERADE show similar performance in key indicators such as the allocation proportion, the number of effectively relieved cases, and their standard deviations, which verifies that they have comparable statistical performance. In addition, compared with Group DBCD, Group ERADE proposed in this paper has a similar average number of allocated people but a smaller standard deviation, demonstrating better stability. Further, compared with the complete randomization design used in the original experiment, both ERADE and Group ERADE significantly increased the number of patients with effective pain relief: the number of patients with a pain reduction of $\ge30\%$ increased from $77$ cases in the original experiment to about $83$ cases, and the number of patients with a pain reduction of $\ge50\%$ increased from $61$ cases to about $66$ cases.

We also explored the robustness of each method under the condition of reducing the sample size. The results show that under delayed and missing responses, Group ERADE proposed in this paper shows similar and satisfactory performance in terms of the allocation proportion, the number of effectively relieved cases, and their standard deviations compared with similar improved ERADE and Group DBCD, further supporting its applicability in practical application scenarios with limited resources or small sample sizes.

\section{Concluding Remarks}\label{sec4}

This paper proposes a new class of group response-adaptive randomization procedures, Group ERADE，based on ERADE framework proposed by Hu, Zhang, and He \cite{hu2009efficient}. A core innovation of this method is that patients are not enrolled one-by-one, but enter the trial in groups according to a preset time period (such as weekly, bi-weekly, or monthly). This design not only retains the advantages of the original ERADE in terms of statistical efficiency and asymptotic optimality but also significantly enhances its practicality in the operation of clinical trials. To further enhance the method's adaptability to real-world scenarios, we propose an improved version of Group ERADE, specifically designed to address the commonly existing problems of delayed response and data missing in clinical trials. Theoretically, we systematically derive the asymptotic properties of the proposed design and its improved version. The results show that under relatively broad regular conditions, Group ERADE has the same asymptotic behavior as the original ERADE, including the strong consistency and asymptotic normality of the allocation proportion. 

In some clinical trials, covariates are important factors that influence treatment allocation and efficacy evaluation. If covariate information is ignored in adaptive designs, it may lead to imbalances between treatments, thereby introducing bias. Taves \cite{taves1974minimization} proposed a minimization allocation method, which effectively reduced the interference of differences in patient characteristics on efficacy comparison by dynamically balancing the distribution of covariates between treatments. On this basis, Pocock and Simon \cite{pocock1975sequential} further developed a sequential allocation strategy for balancing prognostic factors, laying the foundation for the application of covariate information in clinical trial randomization. Atkinson \cite{atkinson1982optimum} explored the biased coin design method in sequential clinical trials with prognostic factors from the perspective of optimal design. Subsequently, Hu and Rosenberger \cite{hu2006theory} expanded the theoretical framework of response-adaptive randomized designs, incorporating both patients' response data and baseline covariate information into the update mechanism of allocation probabilities, thus improving the efficiency and fairness of the design.

From the perspective of group entry, the Group ERADE proposed in this paper has preliminarily verified the feasibility and theoretical properties of the fixed-period group entry method. However, if covariate information is further introduced, several issues remain inadequately resolved, such as how to balance the covariate at the group level, how to use covariates to predict the delayed and missing responses within the group, and how to update the allocation probability including covariates in real time when there are delayed responses.

In summary, future research can be carried out in the following directions: (1) Incorporate covariate information into the allocation probability function of Group ERADE and explore its asymptotic properties within the group entry framework; (2) Develop covariate-based methods for predicting delayed and missing responses and combine them with the group adaptive update mechanism; (3) Expand the existing theory to handle complex situations where covariates change over time or are partially missing. We leave these topics as important content for subsequent research.

\bibliographystyle{elsarticle-num}
\bibliography{ref}

\section{Appendix}

\subsection{Appendix.A}

\begin{table}[h]
	\centering
	\small
	\caption{Results for the binary case with Neyman allocation proportion.}
	\resizebox{\textwidth}{!}{
		\begin{tabular}{cccccccc}
			\hline
			& & \multicolumn{2}{c}{ERADE} & \multicolumn{2}{c}{Group ERADE} & \multicolumn{2}{c}{Group DBCD} \\
			\cline{3-8} 
			$P_1$, $P_2$​ & \makecell[c]{$\nu$\\$\sigma/\sqrt{\tau_j}$} & \makecell[c]{$\mathbb{E}\left({N_{\tau_j,1}}/{\tau_j}\right)$\\(SD)} & \makecell[c]{$\mathbb{E}\left(\text{No.Fail}\right)$\\(SD)} & \makecell[c]{$\mathbb{E}\left({N_{\tau_j,1}}/{\tau_j}\right)$\\(SD)} & \makecell[c]{$\mathbb{E}\left(\text{No.Fail}\right)$\\(SD)} & \makecell[c]{$\mathbb{E}\left({N_{\tau_j,1}}/{\tau_j}\right)$\\(SD)} & \makecell[c]{$\mathbb{E}\left(\text{No.Fail}\right)$\\(SD)}\\\hline
			0.9, 0.7 & \makecell[c]{0.396\\(0.037)} & \makecell[c]{0.391\\(0.048)} & \makecell[c]{0.222\\(0.028)} & \makecell[c]{0.393\\(0.047)} & \makecell[c]{0.221\\(0.028)} & \makecell[c]{0.391\\(0.051)} & \makecell[c]{0.222\\(0.028)} \\\hline
			0.9, 0.5 & \makecell[c]{0.375\\(0.036)} & \makecell[c]{0.369\\(0.047)} & \makecell[c]{0.353\\(0.030)} & \makecell[c]{0.372\\(0.047)} & \makecell[c]{0.350\\(0.030)} & \makecell[c]{0.371\\(0.049)} & \makecell[c]{0.352\\(0.031)} \\\hline
			0.8, 0.8 & \makecell[c]{0.500\\(0.027)} & \makecell[c]{0.500\\(0.031)} & \makecell[c]{0.200\\(0.028)} & \makecell[c]{0.500\\(0.032)} & \makecell[c]{0.200\\(0.029)} & \makecell[c]{0.501\\(0.037)} & \makecell[c]{0.200\\(0.028)} \\\hline
			0.8, 0.6 & \makecell[c]{0.449\\(0.020)} & \makecell[c]{0.449\\(0.023)} & \makecell[c]{0.310\\(0.030)} & \makecell[c]{0.450\\(0.025)} & \makecell[c]{0.311\\(0.031)} & \makecell[c]{0.447\\(0.031)} & \makecell[c]{0.309\\(0.030)} \\\hline
			0.7, 0.5 & \makecell[c]{0.478\\(0.011)} & \makecell[c]{0.478\\(0.013)} & \makecell[c]{0.404\\(0.033)} & \makecell[c]{0.478\\(0.016)} & \makecell[c]{0.404\\(0.032)} & \makecell[c]{0.478\\(0.021)} & \makecell[c]{0.405\\(0.033)} \\\hline
			0.6, 0.4 & \makecell[c]{0.500\\(0.007)} & \makecell[c]{0.500\\(0.009)} & \makecell[c]{0.499\\(0.033)} & \makecell[c]{0.500\\(0.012)} & \makecell[c]{0.500\\(0.033)} & \makecell[c]{0.500\\(0.019)} & \makecell[c]{0.500\\(0.033)} \\\hline
			0.5, 0.5 & \makecell[c]{0.500\\(0.000)} & \makecell[c]{0.500\\(0.006)} & \makecell[c]{0.500\\(0.036)} & \makecell[c]{0.500\\(0.010)} & \makecell[c]{0.500\\(0.036)} & \makecell[c]{0.500\\(0.016)} & \makecell[c]{0.500\\(0.035)} \\\hline
			0.4, 0.3 & \makecell[c]{0.517\\(0.012)} & \makecell[c]{0.517\\(0.014)} & \makecell[c]{0.648\\(0.033)} & \makecell[c]{0.517\\(0.017)} & \makecell[c]{0.649\\(0.033)} & \makecell[c]{0.518\\(0.022)} & \makecell[c]{0.648\\(0.033)} \\\hline
			0.2, 0.2 & \makecell[c]{0.500\\(0.027)} & \makecell[c]{0.500\\(0.030)} & \makecell[c]{0.800\\(0.028)} & \makecell[c]{0.500\\(0.031)} & \makecell[c]{0.799\\(0.028)} & \makecell[c]{0.499\\(0.037)} & \makecell[c]{0.800\\(0.028)} \\\hline
		\end{tabular}
	}
\end{table}

\begin{table}[h]
	\centering
	\small
	\caption{Delayed and missing response: results for the binary case with Neyman target.}
	\resizebox{\textwidth}{!}{
		\begin{tabular}{cccccccc}
			\hline
			& & \multicolumn{2}{c}{ERADE} & \multicolumn{2}{c}{Group ERADE} & \multicolumn{2}{c}{Group DBCD} \\
			\cline{3-8} 
			$\beta$, $P_1$, $P_2$​ & \makecell[c]{$\nu$\\$\sigma/\sqrt{\tau_j}$} & \makecell[c]{$\mathbb{E}\left({N_{\tau_j,1}}/{\tau_j}\right)$\\(SD)} & \makecell[c]{$\mathbb{E}\left(\text{No.Fail}\right)$\\(SD)} & \makecell[c]{$\mathbb{E}\left({N_{\tau_j,1}}/{\tau_j}\right)$\\(SD)} & \makecell[c]{$\mathbb{E}\left(\text{No.Fail}\right)$\\(SD)} & \makecell[c]{$\mathbb{E}\left({N_{\tau_j,1}}/{\tau_j}\right)$\\(SD)} & \makecell[c]{$\mathbb{E}\left(\text{No.Fail}\right)$\\(SD)}\\\hline
			0.0 0.9, 0.7 & \makecell[c]{0.396\\(0.037)} & \makecell[c]{0.389\\(0.054)} & \makecell[c]{0.222\\(0.028)} & \makecell[c]{0.394\\(0.047)} & \makecell[c]{0.221\\(0.028)} & \makecell[c]{0.392\\(0.051)} & \makecell[c]{0.221\\(0.029)} \\\hline
			0.0 0.9, 0.5 & \makecell[c]{0.375\\(0.036)} & \makecell[c]{0.365\\(0.055)} & \makecell[c]{0.354\\(0.033)} & \makecell[c]{0.373\\(0.047)} & \makecell[c]{0.351\\(0.030)} & \makecell[c]{0.371\\(0.049)} & \makecell[c]{0.352\\(0.031)} \\\hline
			0.0 0.8, 0.8 & \makecell[c]{0.500\\(0.027)} & \makecell[c]{0.501\\(0.038)} & \makecell[c]{0.200\\(0.029)} & \makecell[c]{0.501\\(0.032)} & \makecell[c]{0.200\\(0.028)} & \makecell[c]{0.500\\(0.038)} & \makecell[c]{0.200\\(0.028)} \\\hline
			0.0 0.8, 0.6 & \makecell[c]{0.449\\(0.020)} & \makecell[c]{0.447\\(0.030)} & \makecell[c]{0.310\\(0.031)} & \makecell[c]{0.450\\(0.025)} & \makecell[c]{0.311\\(0.030)} & \makecell[c]{0.448\\(0.031)} & \makecell[c]{0.311\\(0.031)} \\\hline
			0.0 0.7, 0.5 & \makecell[c]{0.478\\(0.011)} & \makecell[c]{0.478\\(0.016)} & \makecell[c]{0.405\\(0.032)} & \makecell[c]{0.478\\(0.016)} & \makecell[c]{0.405\\(0.032)} & \makecell[c]{0.478\\(0.022)} & \makecell[c]{0.405\\(0.032)} \\\hline
			0.0 0.6, 0.4 & \makecell[c]{0.500\\(0.007)} & \makecell[c]{0.500\\(0.011)} & \makecell[c]{0.501\\(0.034)} & \makecell[c]{0.500\\(0.013)} & \makecell[c]{0.501\\(0.033)} & \makecell[c]{0.500\\(0.019)} & \makecell[c]{0.500\\(0.033)} \\\hline
			0.0 0.5, 0.5 & \makecell[c]{0.500\\(0.000)} & \makecell[c]{0.500\\(0.007)} & \makecell[c]{0.499\\(0.035)} & \makecell[c]{0.500\\(0.010)} & \makecell[c]{0.499\\(0.035)} & \makecell[c]{0.500\\(0.016)} & \makecell[c]{0.499\\(0.035)} \\\hline
			0.0 0.4, 0.3 & \makecell[c]{0.517\\(0.012)} & \makecell[c]{0.516\\(0.016)} & \makecell[c]{0.649\\(0.033)} & \makecell[c]{0.517\\(0.017)} & \makecell[c]{0.649\\(0.033)} & \makecell[c]{0.518\\(0.022)} & \makecell[c]{0.648\\(0.033)} \\\hline
			0.0 0.2, 0.2 & \makecell[c]{0.500\\(0.027)} & \makecell[c]{0.500\\(0.038)} & \makecell[c]{0.800\\(0.028)} & \makecell[c]{0.499\\(0.033)} & \makecell[c]{0.800\\(0.028)} & \makecell[c]{0.499\\(0.038)} & \makecell[c]{0.800\\(0.029)} \\\hline
			0.1 0.9, 0.7 & \makecell[c]{0.396\\(0.039)} & \makecell[c]{0.388\\(0.059)} & \makecell[c]{0.222\\(0.031)} & \makecell[c]{0.393\\(0.053)} & \makecell[c]{0.221\\(0.030)} & \makecell[c]{0.391\\(0.055)} & \makecell[c]{0.223\\(0.031)} \\\hline
			0.1 0.9, 0.5 & \makecell[c]{0.375\\(0.038)} & \makecell[c]{0.366\\(0.057)} & \makecell[c]{0.353\\(0.034)} & \makecell[c]{0.373\\(0.050)} & \makecell[c]{0.351\\(0.032)} & \makecell[c]{0.369\\(0.054)} & \makecell[c]{0.352\\(0.033)} \\\hline
			0.1 0.8, 0.8 & \makecell[c]{0.500\\(0.028)} & \makecell[c]{0.501\\(0.041)} & \makecell[c]{0.200\\(0.030)} & \makecell[c]{0.500\\(0.037)} & \makecell[c]{0.200\\(0.030)} & \makecell[c]{0.498\\(0.042)} & \makecell[c]{0.200\\(0.029)} \\\hline
			0.1 0.8, 0.6 & \makecell[c]{0.449\\(0.021)} & \makecell[c]{0.447\\(0.033)} & \makecell[c]{0.311\\(0.033)} & \makecell[c]{0.448\\(0.031)} & \makecell[c]{0.310\\(0.032)} & \makecell[c]{0.447\\(0.034)} & \makecell[c]{0.312\\(0.033)} \\\hline
			0.1 0.7, 0.5 & \makecell[c]{0.478\\(0.012)} & \makecell[c]{0.477\\(0.020)} & \makecell[c]{0.403\\(0.034)} & \makecell[c]{0.478\\(0.020)} & \makecell[c]{0.404\\(0.034)} & \makecell[c]{0.477\\(0.026)} & \makecell[c]{0.405\\(0.035)} \\\hline
			0.1 0.6, 0.4 & \makecell[c]{0.500\\(0.008)} & \makecell[c]{0.499\\(0.016)} & \makecell[c]{0.500\\(0.036)} & \makecell[c]{0.500\\(0.017)} & \makecell[c]{0.499\\(0.036)} & \makecell[c]{0.500\\(0.023)} & \makecell[c]{0.499\\(0.036)} \\\hline
			0.1 0.5, 0.5 & \makecell[c]{0.500\\(0.000)} & \makecell[c]{0.500\\(0.013)} & \makecell[c]{0.501\\(0.037)} & \makecell[c]{0.500\\(0.016)} & \makecell[c]{0.500\\(0.037)} & \makecell[c]{0.500\\(0.021)} & \makecell[c]{0.500\\(0.037)} \\\hline
			0.1 0.4, 0.3 & \makecell[c]{0.517\\(0.013)} & \makecell[c]{0.515\\(0.021)} & \makecell[c]{0.648\\(0.036)} & \makecell[c]{0.517\\(0.022)} & \makecell[c]{0.648\\(0.035)} & \makecell[c]{0.517\\(0.026)} & \makecell[c]{0.648\\(0.035)} \\\hline
			0.1 0.2, 0.2 & \makecell[c]{0.500\\(0.028)} & \makecell[c]{0.498\\(0.042)} & \makecell[c]{0.800\\(0.030)} & \makecell[c]{0.501\\(0.037)} & \makecell[c]{0.800\\(0.030)} & \makecell[c]{0.500\\(0.042)} & \makecell[c]{0.799\\(0.030)} \\\hline
			0.2 0.9, 0.7 & \makecell[c]{0.396\\(0.041)} & \makecell[c]{0.388\\(0.062)} & \makecell[c]{0.222\\(0.033)} & \makecell[c]{0.391\\(0.056)} & \makecell[c]{0.221\\(0.032)} & \makecell[c]{0.392\\(0.058)} & \makecell[c]{0.222\\(0.033)} \\\hline
			0.2 0.9, 0.5 & \makecell[c]{0.375\\(0.040)} & \makecell[c]{0.366\\(0.059)} & \makecell[c]{0.354\\(0.037)} & \makecell[c]{0.372\\(0.054)} & \makecell[c]{0.351\\(0.035)} & \makecell[c]{0.371\\(0.055)} & \makecell[c]{0.352\\(0.035)} \\\hline
			0.2 0.8, 0.8 & \makecell[c]{0.500\\(0.030)} & \makecell[c]{0.501\\(0.045)} & \makecell[c]{0.200\\(0.031)} & \makecell[c]{0.500\\(0.043)} & \makecell[c]{0.200\\(0.032)} & \makecell[c]{0.499\\(0.047)} & \makecell[c]{0.200\\(0.032)} \\\hline
			0.2 0.8, 0.6 & \makecell[c]{0.449\\(0.023)} & \makecell[c]{0.446\\(0.038)} & \makecell[c]{0.311\\(0.035)} & \makecell[c]{0.449\\(0.036)} & \makecell[c]{0.310\\(0.035)} & \makecell[c]{0.448\\(0.038)} & \makecell[c]{0.311\\(0.035)} \\\hline
			0.2 0.7, 0.5 & \makecell[c]{0.478\\(0.012)} & \makecell[c]{0.477\\(0.024)} & \makecell[c]{0.404\\(0.037)} & \makecell[c]{0.479\\(0.025)} & \makecell[c]{0.405\\(0.037)} & \makecell[c]{0.477\\(0.030)} & \makecell[c]{0.405\\(0.037)} \\\hline
			0.2 0.6, 0.4 & \makecell[c]{0.500\\(0.008)} & \makecell[c]{0.499\\(0.021)} & \makecell[c]{0.500\\(0.038)} & \makecell[c]{0.500\\(0.023)} & \makecell[c]{0.501\\(0.037)} & \makecell[c]{0.500\\(0.027)} & \makecell[c]{0.500\\(0.037)} \\\hline
			0.2 0.5, 0.5 & \makecell[c]{0.500\\(0.000)} & \makecell[c]{0.500\\(0.019)} & \makecell[c]{0.501\\(0.039)} & \makecell[c]{0.500\\(0.020)} & \makecell[c]{0.499\\(0.039)} & \makecell[c]{0.500\\(0.025)} & \makecell[c]{0.500\\(0.040)} \\\hline
			0.2 0.4, 0.3 & \makecell[c]{0.517\\(0.014)} & \makecell[c]{0.516\\(0.026)} & \makecell[c]{0.648\\(0.037)} & \makecell[c]{0.516\\(0.026)} & \makecell[c]{0.649\\(0.037)} & \makecell[c]{0.518\\(0.030)} & \makecell[c]{0.649\\(0.038)} \\\hline
			0.2 0.2, 0.2 & \makecell[c]{0.500\\(0.030)} & \makecell[c]{0.498\\(0.046)} & \makecell[c]{0.800\\(0.032)} & \makecell[c]{0.500\\(0.043)} & \makecell[c]{0.800\\(0.032)} & \makecell[c]{0.500\\(0.047)} & \makecell[c]{0.800\\(0.032)} \\\hline
		\end{tabular}
	}
\end{table}

\begin{table}[h]
	\centering
	\small
	\caption{Results for the continuous case with Neyman allocation proportion.}
	\resizebox{\textwidth}{!}{
		\begin{tabular}{cccccccc}
			\hline
			& & \multicolumn{2}{c}{ERADE} & \multicolumn{2}{c}{Group ERADE} & \multicolumn{2}{c}{Group DBCD} \\
			\cline{3-8} 
			$\mu_1$,$\mu_2$,$\pi_1$,$\pi_2$​ & \makecell[c]{$\nu$\\$\sigma/\sqrt{\tau_j}$} & \makecell[c]{$\mathbb{E}\left({N_{\tau_j,1}}/{\tau_j}\right)$\\(SD)} & \makecell[c]{$\mathbb{E}\left(\text{Obs.Rsp}\right)$\\(SD)} & \makecell[c]{$\mathbb{E}\left({N_{\tau_j,1}}/{\tau_j}\right)$\\(SD)} & \makecell[c]{$\mathbb{E}\left(\text{Obs.Rsp}\right)$\\(SD)} & \makecell[c]{$\mathbb{E}\left({N_{\tau_j,1}}/{\tau_j}\right)$\\(SD)} & \makecell[c]{$\mathbb{E}\left(\text{Obs.Rsp}\right)$\\(SD)}\\\hline
			13,4,15,2.5 & \makecell[c]{0.615\\(0.024)} & \makecell[c]{0.616\\(0.026)} & \makecell[c]{13.765\\(0.256)} & \makecell[c]{0.613\\(0.027)} & \makecell[c]{13.772\\(0.256)} & \makecell[c]{0.617\\(0.032)} & \makecell[c]{13.757\\(0.259)} \\\hline
			13,2.5,15,4 & \makecell[c]{0.385\\(0.024)} & \makecell[c]{0.385\\(0.026)} & \makecell[c]{14.230\\(0.254)} & \makecell[c]{0.387\\(0.027)} & \makecell[c]{14.226\\(0.251)} & \makecell[c]{0.382\\(0.033)} & \makecell[c]{14.236\\(0.257)} \\\hline
			15,4,17,2.5 & \makecell[c]{0.615\\(0.024)} & \makecell[c]{0.615\\(0.026)} & \makecell[c]{15.772\\(0.248)} & \makecell[c]{0.614\\(0.027)} & \makecell[c]{15.773\\(0.257)} & \makecell[c]{0.617\\(0.032)} & \makecell[c]{15.763\\(0.253)} \\\hline
			15,2.5,17,4 & \makecell[c]{0.385\\(0.024)} & \makecell[c]{0.385\\(0.026)} & \makecell[c]{16.229\\(0.254)} & \makecell[c]{0.386\\(0.027)} & \makecell[c]{16.228\\(0.254)} & \makecell[c]{0.382\\(0.033)} & \makecell[c]{16.240\\(0.254)} \\\hline
		\end{tabular}
	}
\end{table}

\begin{table}[h]
	\centering
	\small
	\caption{Delayed and missing response: results for the continuous case with Neyman allocation proportion.}
	\resizebox{\textwidth}{!}{
		\begin{tabular}{cccccccc}
			\hline
			& & \multicolumn{2}{c}{ERADE} & \multicolumn{2}{c}{Group ERADE} & \multicolumn{2}{c}{Group DBCD} \\
			\cline{3-8} 
			$\beta$,$\mu_1$,$\pi_1$,$\mu_2$,$\pi_2$​ & \makecell[c]{$\nu$\\$\sigma/\sqrt{\tau_j}$} & \makecell[c]{$\mathbb{E}\left({N_{\tau_j,1}}/{\tau_j}\right)$\\(SD)} & \makecell[c]{$\mathbb{E}\left(\text{Obs.Rsp}\right)$\\(SD)} & \makecell[c]{$\mathbb{E}\left({N_{\tau_j,1}}/{\tau_j}\right)$\\(SD)} & \makecell[c]{$\mathbb{E}\left(\text{Obs.Rsp}\right)$\\(SD)} & \makecell[c]{$\mathbb{E}\left({N_{\tau_j,1}}/{\tau_j}\right)$\\(SD)} & \makecell[c]{$\mathbb{E}\left(\text{Obs.Rsp}\right)$\\(SD)}\\\hline
			0.0,13,4,15,2.5 & \makecell[c]{0.615\\(0.024)} & \makecell[c]{0.614\\(0.026)} & \makecell[c]{13.773\\(0.249)} & \makecell[c]{0.613\\(0.028)} & \makecell[c]{13.775\\(0.251)} & \makecell[c]{0.616\\(0.032)} & \makecell[c]{13.766\\(0.258)} \\\hline
			0.0,13,2.5,15,4 & \makecell[c]{0.385\\(0.024)} & \makecell[c]{0.386\\(0.026)} & \makecell[c]{14.226\\(0.250)} & \makecell[c]{0.388\\(0.028)} & \makecell[c]{14.228\\(0.258)} & \makecell[c]{0.385\\(0.033)} & \makecell[c]{14.225\\(0.257)} \\\hline
			0.0,15,4,17,2.5 & \makecell[c]{0.615\\(0.024)} & \makecell[c]{0.614\\(0.026)} & \makecell[c]{15.772\\(0.252)} & \makecell[c]{0.612\\(0.028)} & \makecell[c]{15.775\\(0.254)} & \makecell[c]{0.617\\(0.032)} & \makecell[c]{15.771\\(0.257)} \\\hline
			0.0,15,2.5,17,4 & \makecell[c]{0.385\\(0.024)} & \makecell[c]{0.386\\(0.026)} & \makecell[c]{16.224\\(0.251)} & \makecell[c]{0.388\\(0.027)} & \makecell[c]{16.225\\(0.251)} & \makecell[c]{0.383\\(0.033)} & \makecell[c]{16.236\\(0.252)} \\\hline
			0.2,13,4,15,2.5 & \makecell[c]{0.615\\(0.027)} & \makecell[c]{0.614\\(0.034)} & \makecell[c]{13.767\\(0.284)} & \makecell[c]{0.612\\(0.035)} & \makecell[c]{13.773\\(0.285)} & \makecell[c]{0.617\\(0.040)} & \makecell[c]{13.768\\(0.293)} \\\hline
			0.2,13,2.5,15,4 & \makecell[c]{0.385\\(0.027)} & \makecell[c]{0.386\\(0.034)} & \makecell[c]{14.228\\(0.288)} & \makecell[c]{0.388\\(0.036)} & \makecell[c]{14.225\\(0.282)} & \makecell[c]{0.383\\(0.040)} & \makecell[c]{14.226\\(0.292)} \\\hline
			0.2,15,4,17,2.5 & \makecell[c]{0.615\\(0.027)} & \makecell[c]{0.614\\(0.035)} & \makecell[c]{15.778\\(0.285)} & \makecell[c]{0.612\\(0.036)} & \makecell[c]{15.778\\(0.283)} & \makecell[c]{0.617\\(0.041)} & \makecell[c]{15.766\\(0.288)} \\\hline
			0.2,15,2.5,17,4 & \makecell[c]{0.385\\(0.027)} & \makecell[c]{0.386\\(0.034)} & \makecell[c]{16.227\\(0.285)} & \makecell[c]{0.387\\(0.037)} & \makecell[c]{16.224\\(0.288)} & \makecell[c]{0.384\\(0.041)} & \makecell[c]{16.231\\(0.290)} \\\hline
			0.3,13,4,15,2.5 & \makecell[c]{0.615\\(0.029)} & \makecell[c]{0.614\\(0.040)} & \makecell[c]{13.767\\(0.298)} & \makecell[c]{0.612\\(0.041)} & \makecell[c]{13.780\\(0.310)} & \makecell[c]{0.616\\(0.045)} & \makecell[c]{13.771\\(0.315)} \\\hline
			0.3,13,2.5,15,4 & \makecell[c]{0.385\\(0.029)} & \makecell[c]{0.385\\(0.040)} & \makecell[c]{14.227\\(0.309)} & \makecell[c]{0.388\\(0.040)} & \makecell[c]{14.222\\(0.305)} & \makecell[c]{0.384\\(0.046)} & \makecell[c]{14.235\\(0.309)} \\\hline
			0.3,15,4,17,2.5 & \makecell[c]{0.615\\(0.029)} & \makecell[c]{0.615\\(0.039)} & \makecell[c]{15.772\\(0.307)} & \makecell[c]{0.612\\(0.040)} & \makecell[c]{15.776\\(0.304)} & \makecell[c]{0.617\\(0.046)} & \makecell[c]{15.765\\(0.310)} \\\hline
			0.3,15,2.5,17,4 & \makecell[c]{0.385\\(0.029)} & \makecell[c]{0.386\\(0.039)} & \makecell[c]{16.227\\(0.310)} & \makecell[c]{0.387\\(0.041)} & \makecell[c]{16.220\\(0.305)} & \makecell[c]{0.383\\(0.045)} & \makecell[c]{16.234\\(0.310)} \\\hline
		\end{tabular}
	}
\end{table}

\end{document}